\documentclass[leqno]{amsart}
\usepackage{amssymb,latexsym,amsmath,amscd,epsfig,amsthm, amsfonts}
\theoremstyle{theorem}
\newtheorem{lemma}{Lemma}
\newtheorem{theorem}{Theorem}
\newtheorem{corollary}{Corollary}
\newtheorem*{theorem1}{Theorem 1$'$}

\theoremstyle{remark}
\newtheorem{example}{Example}

\newtheorem*{remark}{Remark}

\newcommand{\R}{{\mathbf R}}
\newcommand{\Z}{{\mathbf Z}}

\newcommand{\const}{{\rm const}}

\newcommand{\dd}{{\rm d}}
\newcommand{\tr}{{\rm tr}}

\newcommand{\supp}{{\rm supp}}

\begin{document}
\author{Eugenia Malinnikova}
\title{Orthonormal sequences in $L^2(\R^d)$ and time frequency localization} 
\address{Department of Mathematical Sciences, Norwegian University of Science and Technology, NO-7491, Trondheim, Norway}
\email{eugenia@math.ntnu.no}

\thanks{The author is supported by the Research
Council of Norway, grants 160192/V30 and  177355/V30.}   
\begin{abstract}
We prove that there does not exist an orthonormal basis $\{b_n\}$ for $L^2(\R)$ such that
the sequences $\{\mu(b_n)\}$, $\{\mu(\widehat{b_n})\}$, and $\{\Delta(b_n)\Delta(\widehat{b_n})\}$ are bounded.
A higher dimensional version of this result that involves generalized dispersions is also obtained. 
The main tool is a time-frequency localization inequality for orthonormal sequences in $L^2(\R^d)$. 
On the other hand, for $d>1$ we construct a basis $\{b_n\}$ for $L^2(\R^d)$ such that the sequences 
$\{\mu(b_n)\}$, $\{\mu(\widehat{b_n})\}$, and $\{\Delta(b_n)\Delta(\widehat{b_n})\}$ are bounded.

\end{abstract}
\subjclass[2000]{42B10, 42C25}
\keywords{Uncertainty principle, Orthonormal basis, Time Frequency localization}
\maketitle

\section{Introduction}

\subsection{Preliminaries and known results}
Let $f\in L^2(\R)$, $\|f\|_2=1$, then
\[
\mu(f)=\int_{\R} t|f(t)|^2dt\quad{\rm and}\quad \Delta(f)=\left(\int_{\R}(t-\mu(f))^2|f|^2dt\right)^{1/2}
\]
are called the time mean of $f$ and the time dispersion of $f$ respectively. The Fourier transform of $f\in L^2(\R)$ is defined by
\[
\widehat{f}(\xi)=\int_{\R}f(x)e^{-2\pi i \xi\cdot x}dx.
\]   
Then $\mu(\widehat{f})$ and $\Delta(\widehat{f})$ are called the frequency mean and frequency dispersion of $f$. 

The classical Heisenberg uncertainty principle reads
\begin{equation}
\label{eq:H2}
\Delta(f)\Delta(\widehat{f})\ge \frac{1}{4\pi}
\end{equation} 
for any $f\in L^2$ with $\|f\|_2=1$.
 
Clearly, \[\|xf(x)\|_2^2=\mu(f)^2+\Delta^2(f).\]The Heisenberg inequality may be also written in the form
\begin{equation}
\label{eq:H1}
\|xf(x)\|_2^2+\|\xi\widehat{f}(\xi)\|_2^2\ge\frac{1}{2\pi}\|f\|^2_2,
\end{equation}
where $f\in L^2(\R)$ is arbitrary.  We refer the reader to survey articles \cite{Ben,FS} 
and monograph \cite{HJ} for various results related to the uncertainty principle.

In this article we consider uncertainty inequalities for orthonormal sequences and bases. 
For some of the first results related to uncertainty inequalities for orthonormal bases we refer the reader to \cite{M} and the references therein.
The construction of Y.~Meyer yields
a wavelet basis $\{\phi_{n}\}_{n=1}^\infty$ for $L^2(\R)$ such that
\[
\sup_{n}\Delta(\phi_{n})\Delta(\widehat{\phi_{n}})<+\infty.\]
A similar basis is obtained for $L^2(\R^d)$ as well, 
see \cite{M} for details. J.~Bourgain proved that there is an orthonormal basis $\{b_n\}_{n=1}^\infty$ for $L^2(\R)$ such that 
\[\Delta(b_n), \Delta(\widehat{b_n})<\frac1{2\sqrt{\pi}}+\epsilon ,\]
see \cite{Bou}. This result was generalized recently by J.~Benedetto and A.~Powell \cite{BP}. The technique was also used by A.~Powell to construct orthonormal bases with other properties, see \cite{P}. 
The result of Bourgain implies that for each $\epsilon>0$ there is an orthonormal basis such that
\[
\sup_n\Delta(b_n)\Delta(\widehat{b_n})<\frac1{4\pi}+\epsilon,\]
so inequality (\ref{eq:H2}) can not be improved for an orthonormal basis.

On the other hand H.~Shapiro proved a number of uncertainty inequalities for orthonormal sequences that are stronger then corresponding inequalities for a single function. For example, using compactness argument, see \cite{Sh}, one can conclude that for any orthonormal sequence $\{f_n\}_{n=1}^\infty$ in $L^2(\R)$
\begin{equation}
\label{eq:Sh1}
\sup_n\|x f_n\|_2^2+\|\xi \widehat{f_n}\|_2^2=+\infty,\end{equation}
so inequality (\ref{eq:H1}) can be refined for an orthonormal sequence. 
It is also proved in \cite{Sh} that if $\phi,\psi\in L^2(\R)$ then any orthonormal sequence $\{f_n\}$ that satisfies
\begin{equation}
\label{eq:Sh2}
|f_n|\le\phi\quad |\widehat{f_n}|\le\psi,
\end{equation}
is finite. This statement is referred to as the Umbrella Theorem.

Quantitative versions of H.~Shapiro's results appeared in a recent article by Ph.~Jaming and A.~Powell,  \cite{JP}, where
in particular the following sharp Mean Dispersion inequality is obtained.\\
\textit{Let $\{e_k\}_{k=0}^n$ be an orthonormal sequence in $L^2(\R)$ then}
\begin{equation}
\label{eq:MD}
\sum_{k=0}^n\left(\mu(e_k)^2+\Delta^2(e_k)+\mu(\widehat{e_k})^2+\Delta^2(\widehat{e_k})\right)\ge \frac{(n+1)^2}{2\pi}.
\end{equation}
The equality is attained for the sequence of Hermite functions, see \cite{JP}. This inequality implies (\ref{eq:Sh1}).
Further, using results of D.~Slepian, H.O.~Pollak, and H.J.~Landau on time-frequency localization, 
Ph.~Jaming and A.~Powell give a quantitative version of the Umbrella Theorem and obtain a number of inequalities for orthonormal basis and also for Riesz basis for $L^2(\R)$.  

%%%%%%%%%%%%%%%%%%%%%%%%%%%

\subsection{Motivation}
Our interest in the uncertainty principles for orthonormal bases started with discussions with Yu.~Lyubarskii and H.~F\"{u}hr that led to the following question:\\
{\slshape{Does there exist an orthonormal basis for $L^2(\R)$ for which both time and frequency means are bounded and the {\textbf{products}} of dispersions are bounded?}}\\
It seems that it is the product of dispersions that has some "physical meaning" in various problems, but we will not speculate on it here.
  
Related results have been obtained by J.~Benedetto in \cite{Ben} and A.~Powell in \cite{P}. It is not difficult to construct an infinite orthonormal sequence with zero time and frequency means and bounded product of  dispersions (see Example 1 in Section \ref{ss:examples}). However the following is true.

\begin{theorem}
\label{th:basis}
There does not exist an orthonormal basis $\{b_n\}_{n=1}^\infty$ for $L^2(\R)$ such that the sequences $\{\mu(b_n)\}_{n=1}^\infty,\ \{\mu(\widehat{b_n})\}_{n=1}^\infty$, and  $\{\Delta(b_n)\Delta(\widehat{b_n})\}_{n=1}^\infty$ are bounded.
\end{theorem}

We remark that another example of a condition on means and dispersions which can be satisfied by an infinite orthonormal sequence but never by an orthonormal  basis was obtained earlier by A.~Powell. It is proved in \cite{P} that there is no orthonormal basis with bounded (both) dispersions and bounded time means. Theorem \ref{th:basis} can be derived from the Mean Dispersion principle. We will not do it, instead we consider a more general problem in higher-dimensional spaces.

%%%%%%%%%%%%%%%%%%%%%%%%%%%%%%%%%%%%%%%%%

\subsection{Main results}
The main goal of this work is to describe a new version of time-frequency localization that yields a number of precise uncertainty inequalities for orthonormal sequences and basis. The results complement those in \cite{J,P,JP}; our approach is simple and works in $\R^d$ for any $d$. We consider the operator that first time-limits the function and then frequency-limits it, following \cite{SP}. However we don't need the theory of Prolate Spheroidal Wave Functions and the celebrated $2WT$ approximation theorem that was used in \cite{JP}. Instead we use an elementary calculation of the trace of the corresponding self-adjoint operator, that can be found for example in \cite{DS, FS}.     
We obtain the following localization inequality.
\begin{theorem}
\label{th:local}
Let $\{\phi_n\}_{n=1}^N$ be an orthonormal system in $L^2(\R^d)$ and let $T$ and $W$ be measurable subsets of $\R^d$. Assume that 
\[
\int_T|\phi_n|^2=1-a^2_n,\quad \int_W|\widehat{\phi_n}|^2=1-b^2_n.
\]
Then 
\[
\sum_{n=1}^N\left(1-\frac 3 2 a_n-\frac 3 2 b_n\right)\le |T||W|.
\]
\end{theorem}

This result provides a quantitative estimate for the Umbrella Theorem in $\R^d$ as well as a number of inequalities for orthonormal sequences. 
For any $p>0$ and $\phi\in L^2(\R^d)$ we define
\[
\tau_p^p(\phi)=\int_{\R^d}|x|^p|\phi(x)|^2dx.\]
Clearly $\tau_p(\phi)\in[0,+\infty]$ and $\tau_p(\phi)>0$ when $\phi\neq 0$. H\"{o}lder's inequality implies that $\tau_p(\phi)\le \tau_q(\phi)$  when $p<q$ and $\|\phi\|_2=1$.  The localization inequality implies the following generalization of the Mean Dispersion principle.

\begin{theorem} 
\label{th:seq}
Let $p$ be positive and let $\{\phi_n\}_{n}$ be an orthonormal sequence in $L^2(\R^d)$ 
then
\begin{equation}
\label{eq:pmd}
\sum_{n=1}^N \left(\tau_p^p(\phi_n)+\tau_p^p(\phi_n)\right)\ge C N^{1+p/2d},
\end{equation}
where $C$ depends on $d$ and $p$ only. Further,
\begin{equation}
\label{eq:seq}
\sum_n\left(\tau_p(\phi_n)+\tau_p(\widehat{\phi_n})\right)^{-2d-\epsilon}<+\infty,
\end{equation}
for any $\epsilon>0$. 
\end{theorem}
In Section \ref{ss:sharp} we show that for $\epsilon=0$ the inequality (\ref{eq:seq}) does not hold in general, we also show that
(\ref{eq:pmd}) is sharp up to a (multiplicative) constant.

As an application of the localization principle we obtain a higher dimensional version of Theorem \ref{th:basis}. We prove the following:
\begin{theorem1}
Let $p>d$ and let  $\{b_n\}_{n=1}^\infty$ be an orthonormal basis for $L^2(\R^d)$. If sequences $\{q_n\}_{n=1}^\infty\subset\R^d$ and $\{r_n\}_{n=1}^\infty\subset \R^d$ are bounded then   
\[
\sup_n\int_{\R^d}|x-q_n|^p|b_n|^2dx\int_{\R^d}|\xi-r_n|^p|\widehat{b_n}|^2d\xi=+\infty.
\]
\end{theorem1}
Clearly Theorem 1 follows from Theorem $1'$. The next theorem shows that the restriction $p>d$ in Theorem $1'$ is necessary.
\begin{theorem}    
\label{th:p<d}
 For $p\le d$ there exists an orthonormal basis $\{b_n\}_{n=1}^\infty$ for $L^2(\R^d)$ and bounded sequences $\{q_n\}_{n=1}^\infty\subset\R^d$ and $\{r_n\}_{n=1}^\infty\subset \R^d$ such that 
\begin{equation}
\label{eq:th4}
\sup_n\int_{\R^d}|x-q_n|^p|b_n|^2dx\int_{\R^d}|\xi-r_n|^p|\widehat{b_n}|^2d\xi<+\infty.
\end{equation}
\end{theorem}
It follows from the prove that for any $a>0$ we may choose such a basis with $q_n=0$ and $r_n\le a$; however we will see below that no basis satisfies 
(\ref{eq:th4}) with $q_n=r_n=0$. 
We use an argument similar to one in  \cite{Bou} to prove the theorem for $p<d$, some additional technical details are needed to make the argument work for $p=d$. The proof of the last theorem implies that \\
\noindent{\it{for $d\ge 2$ there exists a basis for $L^2(\R^d)$ with bounded time and frequency means and bounded products of dispersions.}}
 
 %%%%%%%%%%%%%%%%%%%%%
 
\subsection{Other uncertainty inequalities}
Various versions of the uncertainty principle are known for functions in $L^2(\R^d)$. We consider two particular inequalities.
The first one is a multidimensional version of the inequality of M.G.~Cowling and J.F.~Price. For any $a>0$ there exists $K(a)>0$ such that 
\begin{equation}
\label{eq:dun1}
\||x|^af(x)\|_2\||\xi|^a\widehat{f}(\xi)\|_2\ge K(a)\|f\|_2^2,
\end{equation}
whenever $f\in L^2(\R^d)$, see \cite{Ben}.
The second is a recent inequality of B.~Demange, \cite{D}. Let 
\begin{equation}
\label{eq:weight}
v(x)=|x_1|^{\alpha_1}...\, |x_d|^{\alpha_d},
\end{equation} 
where $\alpha=(\alpha_1,...,\alpha_d)$, and $\alpha_j>0$ for $j=1,...,d$.
There exists $K(\alpha)>0$ such that
\begin{equation}
\label{eq:dun2}
\|v(x) f(x)\|_2\|v(\xi)\widehat{f}(\xi)\|_2\ge K(\alpha)\|f\|_2^2,
\end{equation}
for any $f\in L^2(\R^d)$.

The localization inequality implies that  for 
any orthonormal basis $\{b_n\}_{n=1}^\infty$ for $L^2(\R^d)$
\begin{equation}
\label{eq:dunseq}
\sup_n\||x|^ab_n\|_2\||\xi|^a\widehat{b_n}\|_2=+\infty.
\end{equation}
Remark that the above statement holds for any $a>0$ in contrast to Theorem~$1'$. 
The reason is that we don't allow any time-frequency shifts now, while in Theorem~$1'$ bounded shifts $(q_n,r_n)$ are allowed.  

Inequality (\ref{eq:dunseq}) can be regarded as a version of the uncertainty inequality (\ref{eq:dun1}) for orthonormal bases. 
If we consider (\ref{eq:dun2}) instead of (\ref{eq:dun1}) the situation becomes different.     
We show that for any $v$ of the form (\ref{eq:weight}) there is an orthonormal basis $\{b_n\}_{n=1}^\infty$ for $L^2(\R^d),\ d>1,$ such that
\begin{equation}
\label{eq:hom}
\sup_n\|v(x)b_n\|_2\|v(\xi)\widehat{b_n}\|_2<+\infty.
\end{equation}
Here our argument is a simple version of that of Bourgain, see \cite{Bou}.

The article is organized as follows. Time-frequency localization is discussed in the next Section, we prove Theorem 2 and obtain its various 
applications including Theorem 3; at the end of the section we use the Hermite functions to show that Theorem 3 is sharp. Section 3 is devoted to Theorem $1'$; we  use localization result to show that there is no orthogonal basis with given properties, we also prove (\ref{eq:dunseq}). 
In the last section various orthonormal bases are constructed, we prove Theorem 4 and show that there is a 
basis with bounded means and bounded products of dispersions for $L^2(\R^d)$ when $d>1$; 
finally we construct a basis that satisfies (\ref{eq:hom}).      
      
%%%%%%%%%%%%%%%%%%%%%%%%%%%%%%%%%%%%%%%%%%%%%%%%%%%%%%%%%%%%%%%%%%%%%%%%%%%%%%%%%%%%%%

%%%%%%%%%%%%%%%%%%%%%%%%%%%%%%%%%%%%%%%%%%%%%%%%%%%%%%%%%%%%%%%%%%%%%%%%%%%%%%%%%%%%%%%

\section{Time-frequency localization}

\subsection{Proof of Theorem \ref{th:local}}
Let $T$ and $W$ be two measurable subsets of $\R^d$ and $\{\phi_n\}_{n=1}^N$ be an orthonormal sequence in $L^2(\R^d)$.
Denote by $\chi_T$ the characteristic function of $T$ and
consider the operators $P_T$ and $P_W$ on $L^2(\R^d)$ defined by 
\[P_T(f)=f\chi_{T},\quad  {\rm and}\quad 
P_W(f)(t)=\int_{W}e^{2\pi t\cdot\xi}\widehat{f}(\xi)d\xi=\int_{W}\int_{\R^d}e^{2\pi(t-s)\cdot\xi}f(s)ds d\xi.\]
Then $P_WP_T$ is an integral operator with the kernel (see \cite{DS, FS})
\[
q(s,t)=\chi_T(s)\int_W e^{2\pi i(s-t)\cdot w} dw.\]
A standard calculation in \cite{FS} shows that $P_WP_T$ is a Hilbert-Schmidt operator and  $\|P_WP_T\|^2_{HS}= |W||T|$.
The corresponding self-adjoint operator \[Q=(P_WP_T)^*P_WP_T=P_TP_WP_T\]  is of trace class  (see also \cite{DS}) and   \[\tr(Q)=\|P_WP_T\|_{HS}^2=|W||T|.\] 
Applying Theorem 5.6 from Chapter IV, \cite{GGK}, we obtain
\[
\sum_{n=1}^N\langle Q\phi_n,\phi_n\rangle\le\tr(Q)=|W||T|.
\]
On the other hand,
\[
\langle Q\phi_n,\phi_n\rangle=\langle P_WP_T\phi_n,P_T\phi_n \rangle=\]
\[
\langle \phi_n,\phi_n\rangle-\langle\phi_n-P_T\phi_n,\phi_n\rangle-\langle P_T\phi_n,\phi_n-P_W\phi_n\rangle-\langle P_WP_T\phi_n,\phi_n-P_T\phi\rangle.
\]
Hence
$\langle Q\phi_n,\phi_n\rangle\ge 1-2a_n-b_n$ and
\[
\sum_{n=1}^N(1-2a_n-b_n)\le |W||T|.
\] 
If we consider the operator $\tilde{Q}=(P_TP_W)^*P_TP_W$, we get similarly
\[
\sum_{n=1}^N(1-a_n-2b_n)\le |T||W|.
\]
And the desired time frequency localization inequality follows.

%%%%%%%%%%%%

\subsection{Inequalities for orthonormal sequences} In this section we follow the ideas of \cite{JP}, where various inequalities for orthonormal sequences were derived from a one-dimensional localization principles. We apply the time frequency localization proved in the previous section to obtain rather accurate inequalities.

The following corollary is an immediate consequence of Theorem \ref{th:local}. 
\begin{corollary}
\label{c:concentration}
Let $\{\phi_n\}_{n=1}^N$ be an orthonormal system in $L^2(\R^d)$ such that $\phi_n$ is $\epsilon$-concentrated on a ball $\{|x|<r_0\}$ and 
$\widehat{\phi_n}$ is $\epsilon$-concentrated on a ball $\{|\xi|<\rho_0\}$, for each $n=1,...,N$, i.e.
\[
\int_{|x|<r_0}|\phi_n|^2\ge 1-\epsilon^2,\quad \int_{|\xi|<\rho_0}|\widehat{\phi_n}|^2\ge 1-\epsilon^2.
\]
Then \[ N\le \frac{\pi^dr_0^d\rho_0^d}{(1-3\epsilon)\Gamma\left(\frac d 2+1\right)^2}.\]
\end{corollary}

Another immediate application of the localization inequality is a quantitative version of Shapiro's Umbrella Theorem, we employ localization on arbitrary measurable subsets.
Let $\epsilon$ be positive and $\omega\in L^2(\R^d)$, define 
\[K_\omega(\epsilon)=\inf\left\{|T|:\int_{\R^d\setminus T}|\omega|^2\le\epsilon^2\right\}.\]

\begin{corollary} Let $\phi,\psi\in L^2(\R^d)$ and let $\{f_n\}_{n=1}^N$ be an orthonormal sequence that satisfies
(\ref{eq:Sh2}). Then $N\le(1-3\epsilon)^{-1}K_\phi(\epsilon)K_\psi(\epsilon)$ for each $\epsilon\in (0,\frac1 3)$.
\end{corollary} 
 
\begin{proof}
Let $\omega$ be a non-negative function in $L^2(\R^d)$,we denote by  $\omega^*$  its non-increasing rearrangement defined on $[0,+\infty)$. For each $\epsilon>0$
there exist $T_\omega(\epsilon)$, \[\{\omega>\omega^*(K_\omega(\epsilon))\}\subset T_\omega(\epsilon)\subset \{\omega\ge\omega^*(K_\omega(\epsilon))\},\] such that $|T_\omega(\epsilon)|=K_\omega(\epsilon)$ and
\[
\int_{\R^d\setminus T_\omega(\epsilon)}|\omega|^2=\epsilon^2.\]
 Then for each $n$ we obtain
\[\int_{T_\phi(\epsilon)}|f_n|^2\ge 1-\epsilon^2,\quad
 \int_{T_\psi(\epsilon)}|\widehat{f_n}|^2\ge 1-\epsilon^2.\]
 Thus by Theorem \ref{th:local}, $N(1-3\epsilon)\le K_\phi(\epsilon) K_\psi(\epsilon)$. 
\end{proof}

%%%%%%%%%%%%%%%%%%%%%%%%%%%%%%%%%%%%%%%%%%%%%%%%%%%%%%%%%%%%%%%%%%%%%% 
\subsection{Proof of Theorem \ref{th:seq}}
Let $P_k=\{n: \tau_p(\phi_n)+\tau_p(\widehat{\phi_n})\in[2^{k-1}, 2^k)\}$, where $k$ is an integer. Then 
\[
\int_{\R^d}|x|^p|\phi_n(x)|^2dt\le 2^{kp}\quad{\rm and}\quad
\int_{\R^d}^{\infty}|\xi|^p|\widehat{\phi}(\xi)|^2d\xi\le 2^{kp}
\] whenever $n\in P_k$. It implies that $\phi_n$ is $\frac 1 4$-concentrated on the ball $B(0,2^{k+\frac 2 p})$ both in time and frequency.
The number of elements in $\cup_{j=1}^kP_j$ is less then $c_1(p,d)4^{dk}$, where $c_1(p,d)$ is a constant that does not depend on $k$. There exists integer $k_0$ such that $P_k$ is empty for all $k<k_0$. (The last statement follows also from a theorem of M.~Cowling and J.~Price, see \cite{CP}.)

For given $N>2c_1(p,d)$ choose $k$ such that $2c_1(p,d)4^{dk}>N>2c_1(p,d)4^{d(k-1)}$.  Then at least half of the functions $\{\phi_n\}_{n=1}^N$ does not belong to $\cup_{j=1}^{k-1}P_j$ and we obtain
\[\sum_{n=1}^N(\tau_p^p(\phi_n)+\tau_p^p(\widehat{\phi_n}))\ge \frac{N}{2^{p+1}}2^{kp}\ge a(p,d) N^{1+p/2d}.
\]
For $N<2c_1(p,d)$ we have $\sum_{n=1}^N(\tau_p^p(\phi_n)+\tau_p^p(\widehat{\phi_n}))\ge cN2^{k_0p}$ and (\ref{eq:pmd}) follows.

In order to prove (\ref{eq:seq}) we note that   
\[
 \sum_n\left(\tau_p(\phi_n)+\tau_p(\widehat{\phi_n})\right)^{-2d-\epsilon}\le\sum_{k=k_0}^\infty\sum_{n\in P_k} 2^{(1-k)(2d+\epsilon)}\]
 \[
 \le \sum _{k=k_0}^\infty c(p,d)4^{dk} 2^{(-k)(2d+\epsilon)}<+\infty.
 \]

%%%%%%%%%%%%%%%%%

\subsection{Hermite functions and sharpness of Theorem \ref{th:seq}}\label{ss:sharp}
The Hermite functions are defined by
\[
h_k(t)=\frac{2^{1/4}}{\sqrt{k!}}\left(-\frac {1}{\sqrt{2\pi}}\right)^k 
e^{\pi t^2} \frac{\dd^k}{\dd t^k} e^{-2\pi t^2},\quad k=0,1,2,... .
\]
These functions form an orthonormal basis for $L^2(\R)$ and satisfy $\widehat{h_k}=i^{-k}h_k$,
\[
\mu(h_k)=\mu(\widehat{h_k})=0,\quad \Delta(h_k)=\Delta(\widehat{h_k})=\sqrt{\frac{2k+1}{4\pi}}.\]
It is well known that the Hermite functions are extremal in many problems concerning the uncertainty principle, see for example \cite{FS}, \cite{JP}. 
We will use them to show that (\ref{eq:pmd}) is sharp up to a constant and that inequality (\ref{eq:seq}) does not hold in general when $\epsilon=0$.
Remind that (see, for example, \cite{FS})
\[
xh_k(x)=\frac{\sqrt{k+1}}{2\sqrt{pi}}h_{k+1}(x)+\frac{\sqrt{k}}{2\sqrt{\pi}}h_{k-1}(x).\]
By induction, for positive integer $n$ we obtain
$\tau_{2n}^{2n}(h_k)=\|x^nh-k(x)\|_2^2\le c_nk^n$. Clearly, by H\"{o}lder's inequality, $\tau_p(h_k)=\tau_p(\widehat{h_k})\le \tau_{2n}(h_k)$ for $p<2n$. Thus for each $p>0$ there exist $\kappa_p$ such that
\[
\tau_p(h_k)\le \kappa_p\sqrt{k}.\]
 
We consider the following orthonormal sequence in $L^2(\R^d)$
\[
\phi_I(x)=h_{i_1}(x_1)h_{i_2}(x_2)...h_{i_d}(x_d),\quad {\rm where}\ I=(i_1,...,i_d).
\]
 For $p>0$ we have 
\[\tau_p^p(\phi_I)=\tau_p^p(\widehat{\phi_I})\le c\sum_{m=1}^d \tau_p^p(h_{i_m})\le C\left(\sum_{m=1}^d i_m\right)^{p/2},\]
where $C$ depends on $p$ and $d$ only.

The number of functions in the system $\{\phi_I\}_{|I|\le K}$ is
$N={K+d\choose d}\ge cK^{d}$ and we obtain
\[
\sum_{|I|\le K}(\tau_p^p(\phi_I)+\tau_p^p(\widehat{\pi_I}))\le C \sum_{j=1}^K{j+d-1\choose d-1}j^{p/2}\le C_1K^{d+p/2}\le C_2N^{1+p/2d}.\]
Thus inequality (\ref{eq:pmd}) is sharp up to a multiplicative constant.

Now we look at (\ref{eq:seq}). The sum
$\sum_I\left(\tau_p(\phi_I)+\tau_p(\widehat{\phi_I})\right)^{-a}$
is finite if and only if
\[
+\infty>
\sum_{I}|I|^{-a/2}=\sum_j{j+d-1\choose d-1}j^{-a/2}.\]
This inequality holds if and only if $a>2d$. Thus
\[
\sum_I\left(\tau_p(\phi_I)+\tau_p(\widehat{\phi_I})\right)^{-2d}=+\infty.\]

%%%%%%%%%%%%%%%%%%%%%%%%%%%%%%%%%%%%%%%%%%%%%%%%%%%%%%%%%%%%%%%%%%%%%%%%%%%%%%%%%%%%

\section{Unbounded product of dispersions}

\subsection{Preliminary lemmas}
Our proof of Theorem $1'$  formulated in the Introduction is based on the following lemmas.

\begin{lemma}
\label{l:JK}
Let $p$ be a positive number and $\{\phi_n\}_{n=1}^N$ be an orthonormal system in $L^2(\R^d)$ that satisfies $\tau_p(\phi_n)\le J$ and  $\tau_p(\widehat{\phi_n})\le K$. Then 
\begin{equation}
\label{eq:JK}
N\le c_0(p,d)(JK)^{d}.
\end{equation}
\end{lemma}

\begin{proof}
Clearly, 
each $\phi_n$ is $\epsilon$-concentrated on the ball $\{|x|^p\le \epsilon^{-2}a(p,d)J^p\}$ and  each $\widehat{\phi_n}$ is $\epsilon$-concentrated on the ball $\{|\xi|^p\le \epsilon^{-2}a(p,d)K^p\}$.
Applying Corollary \ref{c:concentration} with $\epsilon=\frac 1 4$, we obtain
$N\le c_0(p,d)(JK)^d.$
\end{proof}

We note  that for $d=1$ and $p=2$ the Mean-Dispersion inequality (\ref{eq:MD}) by Jaming and Powell, see \cite{JP}, implies (\ref{eq:JK}) with $c_0(2,1)=2\pi$. Their results on one-dimensional time frequency localization give also an estimate on $N$ when $d=1$ and $p>0$.

Another lemma we need is known, it follows for example from  Chapter 3.2.5B) in \cite{HJ}, we give a proof of a simple special case here for the convenience of the reader.

\begin{lemma}
\label{l:PW}
Let $b$ and $c$ be  positive numbers, there exists a nonzero function $f$ in $L^2(\R^d)$ such that $f(x)=0$ when $|x|\le b$, and $\widehat{f}(\xi)=0$ when $|\xi|\le c$.
\end{lemma}

\begin{proof} 
It is enough to consider $d=1$, if $g$  is a required function for $d=1$ (and appropriate $b$ and $c$) we take $f(x)=g(x_1)...g(x_n)$. 

Let $PW_c$ be the space of $f\in L^2(\R)$ such that $\widehat{f}(\xi)=0$ when $|\xi|\ge c$. There exists $a$ such that 
\[\|f\|_2\le a\|f\chi_{\{ |x|> b\}}\|_2,\]
for any $f\in PW_c$, see e.g. \cite{K}. It implies that the traces of functions from $PW_c$ on $\{|x|> b\}$ form a closed subspace in $L^2(\{|x|> b\})$ which is obviously not the whole space. Thus there exists $f\in L^2(\{|x|> b\})$ such that
\[
\int_{|x|>b}f(x)\overline{g(x)}dx=0,\]
for any $g\in PW_c$. We extend $f$ by zero on $\{|x|\le b\}$ in order to get the required function. 
\end{proof}

%%%%%%%%%%%%%
\subsection{Proof of Theorem $1'$}
Let $f\in L^2(\R^d)$,  $p>0$, and $a\in\R^d$. We define
\[
\tau_p^p(f,a)=\int_{\R^d}|x-a|^p|f(x)|^2dx.
\]
Assume that $\{b_n\}_{n=1}^\infty$ is an orthonormal basis, and the sequences $\{q_n\}_{n=1}^\infty,\ \{r_n\}_{n=1}^\infty$, and  $\{\tau_p(b_n,q_n)\tau_p(\widehat{b_n},r_n)\}_{n=1}^\infty$ are bounded for some $p>d$. Let 
\[D^2=\sup_n\tau_p(b_n,q_n)\tau_p(\widehat{b_n},r_n),\quad{\rm and}\quad
M=\max\{\sup_n|q_n|,\sup_n|r_n|\}.\]
 We consider 
\[S_{k}=\{b_n: \tau_p(b_n,q_n)\in (D2^{-k},D2^{-k+1}]\},\]
where $k$ is an integer. 
 Clearly, 
 $\{b_n\}_{n=1}^\infty=\cup_{k}S_{k}.$ 
Note  that $\tau_p(\widehat{b_n},r_n)\le 2^{k}D$ for $b_n\in S_{k}$. 
For $b_n\in S_k$ we have 
\[
\tau_p(b_n)\le \tau_p(b_n,q_n)+|q_n|\le 2^{-k+1}D+M
\quad{\rm and}\quad
\tau_p(\widehat{b_n})\le 2^kD+M.\]  
It follows from Lemma \ref{l:JK} that $S_k$ is finite, and if $N_k$ is the number of elements in $S_k$ then
\[N_k\le c_0(p,d)(2^{|k|+1}D+M)^d(D+M)^d\le a(p,d)2^{d|k|}(D+M)^d.\]

%Now let $A_m$ be a subspace of $L^2(\R)$ generated by elements of $S_{j,k}$ and $S'_{j'k'}$ with $k,k'<m$, it is a finite dimensional subspace. 

Let $R$ be a positive number, we take a function $f\in L^2(\R^d),\ \|f\|_2=1$, that  vanishes on $\{|x|<M+R\}$, and whose Fourier transform vanishes on $\{|\xi|<M+R\}$, see   Lemma \ref{l:PW}. Then we have 
\begin{equation}
\label{eq:pars}
1=\|f\|^2= \sum_{k}\sum_{b_n\in S_{k}}|\langle f,b_n\rangle|^2.
\end{equation}
Now if $b_n\in S_k$, $k>0$ then 
\[
\left|\langle f,b_n\rangle\right|\le\int_{|x|>M+R}|f(x)||b_n(x)|dx\le R^{-p/2}\int_{-\infty}^\infty|x-q_n|^{p/2}|f(x)||b_n(x)|dx\]
\begin{equation}
\label{eq:coef1}
\le  R^{-p/2}\|f\|_2 \tau_p^{p/2}(b_n,q_n)
\le \left(\frac{D}{2^{k-1}R}\right)^{p/2}.
\end{equation}
Similarly, for  $b_n'\in S_{k'}$, $k'\le 0$, we have
\begin{equation}
\label{eq:coef2}
\left|\langle f,b_{n'}\rangle\right|=\left|\langle \widehat{f},\widehat {b_{n'}}\rangle\right|\le \left(\frac{2^{k'}D}{R}\right)^{p/2}.
\end{equation}

 Combining the inequalities (\ref{eq:pars}), (\ref{eq:coef1}), and (\ref{eq:coef2}), we get
 \[
1\le\sum_{k=1}^{\infty}\left(\frac{D}{2^{k-1}R}\right)^pN_k+\sum_{k=0}^{\infty}\left(\frac{D}{2^kR}\right)^pN_{-k}\le
\frac{A(p,d,D,M)}{R^p}\sum_{k=0}^\infty 2^{k(d-p)} .\]
Choosing $R$  large enough, we get a contradiction. The theorem is proved. 

%%%%%%%%%%%%%%%%%%%
\subsection{Another unbounded product}
We complete this section by proving inequality (\ref{eq:dunseq}).
 
\begin{theorem}
If  $\{b_n\}_{n=1}^\infty$ is an orthonormal basis for $L^2(\R^d)$ and $p$ is positive then 
\[
\sup_n \tau_p(b_n)\tau_p(\widehat{b_n})=\infty.\] 

\end{theorem}

\begin{proof} Assume that there exists an orthonormal basis such that $\tau_p(b_n)\tau_p(\widehat{b_n})\le C^2$.
Let  
\[A_k=\{b_n: \tau_p(b_n)\in (2^{-k}C, 2^{-k+1}C]\},\] 
where $k$ is integer.
Clearly for $b_n\in A_k$ we have $\tau_p(\widehat{b_n})\le C2^k$. Then each $b_n\in A_k$ is $\frac 1 4$-concentrated on the ball
 $\{|x|^p<c_2(p,d)C2^{-k}\}$ and $\widehat{b_n}$ is $\frac 1 4$-concentrated on the all $\{|\xi|^p<c_2(p,d)C2^k\}$.
Thus the number of elements in $A_k$ is bounded by a constant that does not depend on $k$. 
Let once again use Lemma \ref{l:PW}, we take a function $f$ in $L^2(\R^d)$ that vanishes on $B(0,R)$ with its Fourier transform.
When $k\ge 0$  and $b_n\in A_k$ we get
\[\langle f,b_n\rangle^2\le R^{-p}\tau_p^p(b_n)\le 2CR^{-p}2^{-kp}.\]
When $k<0$ and $b_n\in A_k$ similarly
\[\langle f, b_n\rangle^2=\langle\widehat{f},\widehat{b_n}\rangle^2\le CR^{-p}2^{kp}.\]
We complete the proof by arguments similar to ones used in the previous theorem.
\end{proof}

%%%%%%%%%%%%%%%%%%%%%%%%%%%%%%%%%%%%%%%%%%%%%%%%%%%%%%%%%%%%%%%%%%%%%%%%%%%%%%%%%%%%%%%%%%%%%%%%%%%%%%%%%%%%%%%%%%%%%%%%%%%%%%%%%

\section{Existence of some orthonormal bases for $L^2(\R^d)$}

\subsection{Orthonormal sequences in one dimension}\label{ss:examples}
We start with two examples of orthonormal sequences in $L^2(\R)$. 
\begin{example}
Let $\phi$ be a real-valued even  $C^\infty$-function, $\supp(\phi)\subset{[-2,-1]\cup[1,2]}$ and $\|\phi\|_2=1$. Then $\Delta(\phi),\Delta(\widehat{\phi})<+\infty$. Consider 
$\phi_n(x)=2^{n/2}\phi(2^nx)$, where $n$ is integer, then $\{\phi_n\}_n$ form an orthonormal sequence such that
\[
\mu(\phi_n)=\mu(\widehat{\phi_n})=0,\quad \Delta(\phi_n)\Delta(\widehat{\phi_n})=c.
\] 
This is an example of an infinite orthonormal sequence with zero means and bounded product of dispersions.
\end{example}

There is also an example of an orthonormal \textit{basis} for $L^2(\R)$ with bounded frequency means and bounded product of dispersions.
\begin{example}
There exists a real function $\psi$ and a corresponding wavelet basis  $\psi_{m,n}(t)=2^{m/2}\psi(2^{m}t-n)$; such that $\Delta(\psi)<+\infty$ and $\Delta(\widehat{\psi})<+\infty$, 
see \cite{M}, \cite{Ben}. One has
\[\mu(\psi_{m,n})=2^{-m}(\mu(\psi)+n);\quad \mu(\widehat{\psi_{m,n}})=2^{m}\mu(\widehat{\psi}).\]
\[\Delta(\psi_{m,n})=2^{-m}\Delta(\psi);\quad \Delta(\widehat{\psi}_{m,n})=2^{m}\Delta(\widehat{\psi}).\]
Thus $\Delta(\psi_{m,n})\Delta(\widehat{\psi_{m,n}})=c$ and $\mu(\widehat{\psi_{m,n}})=0$, since $\widehat{\psi}(-\xi)=\overline{\widehat{\psi}(\xi)}$. However the sequence $\mu(\psi_{m,n})$ is unbounded.
%(Interesting estimates of $c$ for the Meyer wavelets are obtained in \cite{Bat,L}.) 
\end{example}

%%%%%%%%%%%%%%%%%%%%%%%%%%%%%%%%%%%%%%%%%%%%%%%%%%%%%%

\subsection{Some orthonormal sequences in higher dimensions}
\label{ss:ortseq} In this section we obtain preliminary results that we use later to prove Theorem \ref{th:p<d}. First
we construct an orthonormal sequence with required properties that is large in some sense.

Let $\chi$ be the characteristic function of the cube 
\[\{x=(x_1,...,x_d):\ 5/2<x_m<7/2, m=1,...,d\}\]
 and $\omega$ be a smooth radial function supported in $B(0,1/2)$. Then $\phi=\chi*\omega$ is a smooth non-negative function, define $\psi(x)=\phi(x)^{1/2}$. Then 
 \[
\langle \psi(x),e^{2\pi i x\cdot b}\psi(x)\rangle=\widehat{\phi}(2\pi b)=\widehat{\chi}(2\pi b)\widehat{\omega}(2\pi b)=0
\]
whenever $b\in \Z^d,\ b\neq 0$. Further define 
$\Psi(x)=a\psi\left(\frac{x}{2}\right)$, where $a$ is chosen
such that $\|\Psi\|_2=1$. Clearly, \[\supp (\Psi)\subset \{x=(x_1,...,x_d): 1/2<x_m<1\}.\] We have 
\begin{equation}
\label{eq:coef0}
\langle \Psi(x),e^{2\pi i x\cdot b}\Psi(x)\rangle=0,
\end{equation}
when $2b\in \Z^d,\ b\neq 0$.
For every positive integer $s$  and every $j=(j_1,...,j_d)\in \Z^d$, $|j_l|\le 2^s$, we define $\Psi_{j,s}$ by 
\begin{equation}
\label{eq:sjpsi}
\Psi_{j,s}(x)=2^{-ds/2}e^{2\pi i j\cdot 2^{-s}x}\Psi(2^{-s}x).
\end{equation}

\begin{lemma}
\label{l:psi}
 Let $\Psi_{j,s}$ be defined as above. Then 
\begin{equation}
\label{eq:supp}
\supp(\Psi_{j,s})\subset\{x=(x_1,..., x_d): 2^{s-1}<x_m<2^s, m=1,...,d\},
\end{equation}
the sequence $\{\Psi_{j,s}\}_{j,s}$ is orthonormal, and for each $p>0$
there exist $C_1, C_2>0$ such that 
\begin{equation}
\label{eq:tauPsi} 
\tau_p(\Psi_{j,s})=2^sC_1,\quad\tau_p(\widehat{\Psi_{j,s}},2^{-s}j)=2^{-s}C_2 .
\end{equation}
\end{lemma} 

\begin{proof}
The supports of $\Psi_{j,s}$ and $\Psi_{j',s'}$ are disjoint when $s\neq s'$. When $s=s'$ and $j\neq j'$ we have by (\ref{eq:coef0})
\[
\langle \Psi_{j,s},\Psi_{j',s}\rangle=2^{-ds}\langle \Psi(2^{-s}x), e^{2\pi (j'-j)\cdot 2^{-s}x}\Psi(2^{-s}x)\rangle=0.
\]
Further, for any $p>0$,
\[
\tau_p^p(\Psi_{j,s})=\int_{\R^d}|x|^p2^{-ds}|\Psi(2^{-s}x)|^2dx=2^{sp}\tau_p^p(\Psi)=2^{sp}C_1^p.
\]
Clearly $\widehat{\Psi_{j,s}}=2^{ds/2}\widehat{\Psi}(2^s\xi-j)$ and
\[
\tau_p^p(\widehat{\Psi_{j,s}},2^{-s}j)=\int_{\R^d}|\xi-2^{-s}j|^p|\widehat{\Psi_{j,s}}(\xi)|^2d\xi=2^{-sp}\int_{\R^d}|\eta|^p|\widehat{\Psi}(\eta)|^2d\eta=2^{-sp}\tau_p^p(\widehat{\Psi}).
\]
\end{proof}

\begin{remark}
We enumerate $j$ for fixed $s$ as $\{j(n,s)\}_{n=1}^{J_s}$, where $J_s=(2^{s+1}+1)^d>2^{sd}$ and write $\Psi_{n,s}=\Psi_{j(n,s),s}$ for $n=1,...,J_s$.
\end{remark}

\begin{lemma}
\label{l:fourierder}
Let $\Psi_{j,s}$ satisfy (\ref{eq:sjpsi}), where $\supp{\Psi}\subset [-1,1]^d$ and let $q$ be a positive integer. Then there exists $A(q,\Psi)$ such that for any $s$ and  any $R(x)=\sum_j\alpha_j\Psi_{j,s}$ the following inequality holds 
\[
\|\partial_m^qR\|_2^2\le A(\Psi,q)\sum_j|\alpha_j|^2,
\]
here $m\in\{1,...,d\}$ and $\partial_m$ denotes the partial derivative with respect to $x_m$.
\end{lemma}
\begin{proof}
By (\ref{eq:sjpsi}) we have $R(x)=2^{-sd/2}P(2^{-s}x)\Psi(2^{-s}x)$, where 
$P(y)=\sum_j\alpha_j e^{2\pi i j\cdot y}$ is a trigonometric polynomial. We have
\begin{equation*}
%\label{eq:deriv}
\|\partial_m^qR\|_2^2=2^{-2sq}\|\partial_m^q(P\Psi)\|_2^2\le 2^{-2sq}A(q)\sum_{r=0}^q\|\partial_m^{r}P\partial_m^{q-r}\Psi\|_2^2. 
\end{equation*}
Now  $\supp(\Psi)\subset Q=[-1,1]^d$, the functions $e^{2\pi i j\cdot y}$ are orthogonal on this cube and have the same norms. We obtain
\begin{equation*}
%\label{eq:orth} 
 \|\partial_m^rP\partial_t^{q-r}\Psi\|_2^2\le A_1\|\chi_{Q}\partial_m^rP\|_2^2\le A_2\sum_j|\alpha_j(2\pi j_m)^r|^2\le A_32^{2sr}\sum_{j}|\alpha_j|^2,
\end{equation*}
where $A_1,A_2,A_3$ depend on $q$ and $\Psi$. The required inequality follows. 
\end{proof}

%%%%%%%%%%%%%%%%%%%%%%%%%%%%%%%%%%%%%%%%%%%%%%%%%%%%%%%%%%%%%

\subsection{Proof of Theorem \ref{th:p<d}. Case I: $p<d$}
We use the  construction described in \cite{Bou} to replace an orthonormal sequence by a basis, we repeat the details of the construction for the convenience of the reader. 

Let a sequence $\{f_k\}_{k=1}^\infty$ of smooth functions with compact supports be dense on the unit sphere in $L^2(\R^d)$. The basis $\{b_n\}_{n=1}^\infty$ is obtained as $\cup_k \mathcal{B}_k$, where each  $\mathcal{B}_k$ is a finite orthonormal system.
Suppose that a finite  orthonormal system of smooth functions $\mathcal{B}_1,...,\mathcal{B}_{k-1}$ with compact supports is obtained, let $B_{k-1}$ be the linear span of these functions, we put $B_{0}=\{0\}$. We define
\[
f=f_k-P_{B_{k-1}}f_k.\]
Then $f$ is a smooth function with compact support and $I_p(f), I_p(\widehat{f})<+\infty$, where
\begin{equation}
\label{eq:Ip}
\int_{\R^d}(|x|+1)^p|f(x)|^2 dx=I_p(f)<\infty.
\end{equation}

We take $s$ big enough so that the supports of $\Psi_{j,s}$ for all $j=(j_1,..j_d)\in \Z^d, |j_l|\le 2^s,$ do not intersect the supports of $f$ and of functions in $B_{k-1}$ and 
\begin{equation}
\label{eq:choices}
J_s\ge2^{ds}>\max\{I_p(f),2^{sp}I_p(\widehat{f})\},
\end{equation} 
here we use the condition $p<d$. We enumerate $\Psi_{j,s}$ as remarked above. Following \cite{Bou} further, define
\begin{eqnarray}
\label{eq:beta}
\beta_1&=&\frac{\theta}{\sqrt{J_s}}f+\gamma_1\Psi_{1,s},   \notag \\  
\beta_2&=&\frac{\theta}{\sqrt{J_s}}f+\sigma_1\Psi_{1,s}+\gamma_2\Psi_{2,s}, \notag\\
\beta_3&=&\frac{\theta}{\sqrt{J_s}}f+\sigma_1\Psi_{1,s}+\sigma_2\Psi_{2,s}+\gamma_3\Psi_{3,s},\\
...\notag \\
\beta_{J_s}&=&\frac{\theta}{\sqrt{J_s}}f+\sigma_1\Psi_{1,s}+...+\sigma_{J_s-1}\Psi_{J_s-1,s}+\gamma_{J_s}\Psi_{J_s,s}, \notag
\end{eqnarray}
here $\theta\in(0,1/2)$.
Clearly, $\beta_l$ are orthogonal to $B_{k-1}$.
The constants $\sigma_1,...,\sigma_{J_s-1}$ and $\gamma_1,...,\gamma_{J_s}$ are chosen to make $\{\beta_l\}_{l=1}^{J_s}$ an orthonormal sequence.
Thus
\begin{equation}
\label{eq:gamma_sigma}
\gamma_l^2=1-\frac{\theta^2}{J_s}\|f\|_2^2-\sum_{n=1}^{l-1}\sigma_n^2,\quad
\sigma_l\gamma_l=-\frac{\theta^2}{J_s}\|f\|_2^2-\sum_{n=1}^{l-1}\sigma_n^2.
\end{equation}
Clearly, $|\gamma_l|<1$ and, by induction, one has
\begin{equation}
\label{eq:gsest}
|\gamma_l| \ge 1-\frac{2\theta^2}{J_s}\quad{\rm and}\quad |\sigma_l|\le \frac{\theta}{J_s}.
 \end{equation}
We take $\theta=1/4$ and 
 estimate $\tau_p(\beta_l)$ first
\begin{equation*}
\tau_p^p(\beta_l) 
 \le \frac{3}{16J_s}\int_{\R^d}|x|^p|f|^2dx+3\int_{\R^d}|x|^p|\Psi_{l,s}|^2dx+
 3\int_{\R^d}|x|^p\left|\sum_{n=1}^{l-1}\sigma_n\Psi_{n,s}\right|^2dx.
 \end{equation*}
 The first term is bounded by $I_p(f)J_s^{-1}<1<2^{sp}$, the second term is less than $2^{sp}\const$, see (\ref{eq:tauPsi}). To estimate the last term, note that $\supp(\Psi) \in [-1,1]^d$ and (\ref{eq:sjpsi}) implies
\[
\int_{\R^d}|x|^p \left|\sum_{n=1}^{l-1}\sigma_n\Psi_{n,s}\right|^2dx=2^{sp}\int_{\R^d}|y|^p\left|\sum_{n=1}^{l-1}\sigma_ne^{2\pi i j(n,s)\cdot y}\right|^2|\Psi(y)|^2dy\]
\[\le
2^{sp}d^{p/2}\int_{\R^d}\left|\sum_{n=1}^{l-1}\sigma_ne^{2\pi i j(n,s)\cdot y}\right|^2|\Psi(y)|^2dy=2^{sp}d^{p/2}\sum_{n=1}^{l-1}|\sigma_n|^2<2^{sp}d^{p/2}.\]
Thus $\tau_p(\beta_l)\le C2^{s}$, where $C$ depends on $p, d$, and on $\Psi$.

For the Fourier transform we estimate $\tau_p(\widehat{\beta_l},2^{-s}j_l)$, where $j_l=j(l,s)$ was defined in Remark in Section \ref{ss:ortseq}. We have
 \[
 \tau_p^p(\widehat{\beta_l},2^{-s}j_l)=\int_{\R^d}|\xi-2^{-s}j_l|^p|\widehat{\beta_l}|^2d\xi\le \frac{3}{16J_s}\int_{\R^d}(|\xi|+1)^p|\widehat{f}|^2d\xi+\]
 \begin{equation}
 \label{eq:three1}
 3\int_{\R_d}|\xi-2^{-s}j_l|^p|\widehat{\Psi_{l,s}}|^2d\xi+
 3\int_{\R^d}|\xi-2^{-s}j_l|^p\left|\sum_{n=1}^{l-1}\sigma_n\widehat{\Psi_{n,s}}\right|^2d\xi.
 \end{equation}
The first term is bounded by $I_p(\widehat{f})J_s^{-1}$ and is less than $2^{-sp}$ due to our choice of $s$. Inequality (\ref{eq:tauPsi}) implies that the second term is less than $2^{-sp}C_2$. We want to show that the third term is small enough. We have
\begin{equation}
\label{eq:2d}
\int_{\R^d}|\xi-2^{-s}j_l|^p\left|\sum_{n=1}^{l-1}\sigma_n\widehat{\Psi_{n,s}}\right|^2d\xi\le A(d)\int_{\R^d}\left(1+\sum_{m=1}^d\xi_m^{2d}\right)|\widehat{R_l}|^2d\xi,
\end{equation}
where 
$R_l(x)=\sum_{n=1}^{l-1}\sigma_n\Psi_{n,s}(x).$
Clearly, 
\begin{equation}
\label{eq:norm}
\|\widehat{R_l}\|_2^2=\|R_l\|_2^2=\sum_{n=1}^{l-1}|\sigma_n|^2\le J_s^{-1}<2^{-sd}.
\end{equation} 
Further,
\begin{equation}
\label{eq:fourier}
\sum_{m=1}^d \int_{\R^d}\xi_m^{2d}|\widehat{R_l}|^2d\xi=(2\pi)^{-2d}\sum_{m=1}^d\|\partial_m^d R_l\|_2^2.
\end{equation}
Using Lemma \ref{l:fourierder}, we obtain 
\begin{equation}
\label{eq:partial}
\|\partial_m^d R_l\|_2^2\le A(\Psi, d)\sum_{n=1}^{l-1}|\sigma_n|^2\le A(\Psi,d)2^{-sd}
\end{equation}

Finally, combining (\ref{eq:three1}-\ref{eq:partial}), we get
\[
\tau_p^p(\widehat{\beta_l}, 2^{-s}j)\le C(\Psi,d)2^{-ps}.\]
Thus $\tau_p(\beta_l)\tau_p(\widehat{\beta_l},2^{-s}j)< C$. We set $\mathcal{B}_k=\{\beta_1,...,\beta_{J_s}\}$ and continue the procedure. 

We want to check that the resulting orthonormal sequence is complete, once again we follow \cite{Bou}. First, 
\[P_{B_k}f_k=P_{B_{k-1}}f_k+\sum_{l=1}^{J_s}\langle f,\beta_l\rangle \beta_l\] and
\[\|P_{B_k}f_k\|^2\ge\|f_k-f\|_2^2+\frac1{16}\|f\|_2^4\ge\frac1{16}.\]
Suppose that the orthonormal sequence $\cup_{k=1}^\infty\mathcal{B}_k$ is not complete; let $B$ be its closed span. Then there exists $g\in L^2(\R^d)$ such that $\|g\|_2=1$ and $g$ is orthogonal to $B$. For some $k$ we have $\|g-f_k\|_2<1/4$ since $\{f_k\}$ is a dense sequence on the unit sphere of $L^2(\R^d)$. Then we obtain a contradiction
\[
\frac1{16}\le\|P_{B_{k}}(f_k)\|_2^2\le\|P_Bf_k\|_2^2=\|P_B(f_k-g)\|_2^2\le\|f_k-g\|_2^2<\frac 1{16}.\] 

%%%%%%%%%%%%%%%%%%%%%%%%%%%%%%%%%%%%%%%%%%%%%%%%%%%%%%%%%%%%%%%%%%%%%%%%%%%%

\subsection{Proof of Theorem \ref{th:p<d}. Case II: $p=d$}
Our argument in the preceding section does not work when $p=d$. We use strict inequality when on each step we take a function $f_k$, construct $f=f_k-P_{B_{k-1}}f_k$ and choose $s$ that satisfies (\ref{eq:choices}). We will improve the argument to obtain the result for $p=d$. 

As before, we start with a sequence $\{f_k\}_{k=1}^\infty$ of smooth functions with compact supports that is dense on the unit sphere in $L^2(\R^d)$. The basis $\{b_n\}$ is obtained as $\cup_k \mathcal{B}_k$, where  $\mathcal{B}_k$ is a finite orthonormal system; but this time we will have \[\mathcal{B}_k=\cup_{t=1}^{T_k}\mathcal{C}_{t,k}.\]
Suppose that an orthonormal system of functions $\mathcal{B}_1,...,\mathcal{B}_{k-1},\mathcal{C}_{1,k},...,\mathcal{C}_{t-1,k}$ with compact supports is obtained, let $B_{t,k-1}$ be the linear span of these functions. We define
\[
g_t=f_k-P_{B_{t,k-1}}f_k.\]
As before, $g_t$ is a smooth function with compact support and using notation (\ref{eq:Ip}) we have $I_d(g_t), I_d(\widehat{g_t})<+\infty$.
We define also 
\[I(g_t)=\sum_{m=1}^d\int_{\R^d}|\xi_m|^{2d}|\widehat{g_t}(\xi)|^2 d\xi,\] then we have 
$I_d(\widehat{g_t})\le c(d)(1+I(g_t))$.

We take $s$ big enough so that the supports of $\Psi_{j,s}$ do not intersect the supports of $g_t$ and of functions in $B_{t,k-1}$ and also 
\[J_s\ge 2^{ds}>I_d(g_t).\]
Further we repeat construction (\ref{eq:beta}-\ref{eq:gsest}) with $f=g_t$ and
\[\theta=\theta_t=(4+I(g_t))^{-1/2}\in(0,1/2).\] We set $\mathcal{C}_{t,k}=\{\beta_1,...,\beta_{J_s}\}.$

As earlier $\beta_l$ are orthogonal to $B_{t,k-1}$. We estimate $\tau_d(\beta_l)$ as above, $\tau_d(\beta_l)\le C2^{s}$.
For the Fourier transform we once again estimate $\tau_d(\widehat{\beta_l},2^{-s}j_l)$. We have as before
 \[
 \tau_d^d(\widehat{\beta_l},2^{-s}j_l)=\int_{\R^d}|\xi-2^{-s}j_l|^d|\widehat{\beta_l}|^2d\xi\le \frac{3\theta^2}{J_s}\int_{\R^d}(|\xi|+1)^d|\widehat{g_t}|^2d\xi\]
 \begin{equation*}
 +3\int_{\R_d}|\xi-2^{-s}j_l|^d|\widehat{\Psi_{l,s}}|^2d\xi+
 3\int_{\R^d}|\xi-2^{-s}j_l|^d|\sum_{n=1}^{l-1}\sigma_n\widehat{\Psi_{n,s}}|^2d\xi.
 \end{equation*}
The first term is bounded by $3\theta_t^2I_d(\widehat{g_t})J_s^{-1}$. Using our choice of $\theta_t$ and repeating 
estimates for the second and third terms from Case I, we get
\[
\tau_d(\widehat{\beta_l}, 2^{-s}j)\le C(\Psi)2^{-s}.\]
Thus $\tau_d(\beta_l)\tau_d(\widehat{\beta_l},2^{-s}j)< C$. 

Now, in contrast to Case I, the projection of $f_k$ onto the subspace spanned by $\mathcal{B}_1,...,\mathcal{B}_{k-1},\mathcal{C}_{1,k},...,\mathcal{C}_{t-1,k},\mathcal{C}_{t,k}$ could be small if $\theta_t$ is small. So we use the same function $f_k$ again to continue the procedure. Note that
\[
g_{t+1}=f_k-P_{B_{t+1,k-1}}f_k=g_t-\sum_{n=1}^{J_s}\langle g_t,\beta_n\rangle \beta_n=
g_t-\frac{\theta_t}{\sqrt{J_s}}\|g_t\|_2^2\sum_{n=1}^{J_s}\beta_n.\]
And, inserting formulas for $\beta_n$, we obtain
\begin{equation*}
%\label{eq:ft}
g_{t+1}=g_t(1-\theta_t^2\|g_t\|_2^2)-
\sum_{n=1}^{J_s}\frac{\theta_t}{\sqrt{J_s}}\|g_t\|_2^2(\gamma_n+(J_s-n)\sigma_n)\Psi_{n,s}
=\lambda_tg_t+\sum_{n=1}^{J_s}\kappa_n\Psi_{n,s},
\end{equation*}
where $\lambda_t=1-\theta^2\|g_t\|_2^2<1$ and 
\begin{equation}
\label{eq:kappa}
|\kappa_n|=\frac{\theta_t}{\sqrt{J_s}}\|g_t\|_2^2|\gamma_n+(J_s-n)\sigma_n|<\frac{2\theta_t}{\sqrt{J_s}}\|g_t\|_2^2.
\end{equation}
We have
\begin{equation*}
%\label{eq:Ift}
I(g_{t+1})=\sum_{m=1}^d\int_{\R^d}|\xi_m|^{2d}|\widehat{g_{t+1}}|^2d\xi=(2\pi)^{-2d}
\sum_{m=1}^d\|\partial_m^dg_{t+1}\|_2^2.
\end{equation*}
According to our choice of $s$ the support of $g_t$ does not intersect the supports of $\Psi_{n,s}$. Then 
\[
I(g_{t+1})=|\lambda_t|^2I(g_t)+
(2\pi)^{-2d}\sum_{m=1}^d\|\partial_m^d(\sum_{n=1}^{J_s}\kappa_n\Psi_{n,s})\|_2^2. 
\]  
Now Lemma \ref{l:fourierder} implies
\[
\|\partial_m^d(\sum_{n=1}^{J_s}\kappa_n\Psi_{n,s})\|_2^2\le A(d,\Psi)\sum_{n=1}^{J_s}|\kappa_n|^2.
\] 
Combining the last two inequalities with (\ref{eq:kappa}) we obtain
\[
I(g_{t+1})\le\lambda_t^2I(g_t)+B(d,\Psi)\theta_t^2\|g_t\|_2^4\le I(g_1)+B(d,\Psi)\sum_{u=1}^t\theta_u^2\|g_u\|_2^4.
\]
We note also that 
\[
1\ge \|P_{B_{t,k-1}}f_k\|_2^2=\|f_k-g_{t+1}\|^2=\|f_k-g_t\|_2^2+\theta_t^2\|g_t\|_2^4=
\|f_k-g_1\|_2^2+\sum_{u=1}^t\theta_u^2\|g_u\|_2^4.\] 
Thus $I(g_{t+1})\le I(g_1)+B(d,\Psi)$ and $\theta_t^2>(4+I(g_1)+B(d,\Psi))^{-1}$ for each $t$. 
We take $T_k>4+I(g_1)+B(d,\Psi)$ and for $B_k=B_{T_k,k-1}$ we get
\[
\|P_{B_k}f_k\|_2^2=\|f_k-g_1\|_2^2+\sum_{u=1}^{T_k}\theta_u^2\|g_u\|_2^4.
\]
If for each $u$ we have $\|g_u\|_2>1/2$, then 
\[
\|P_{B_k}f_k\|_2^2\ge \frac{T_k}{16}(4+I(g_1)+B(d,\Psi))^{-1}>\frac{1}{16}.
\]
If $\|g_u\|_2<1/2$ for some $u$, then
\[
\|P_{B_k}f_k\|_2^2\ge \|f_k-g_u\|_2^2\ge\frac{1}{4}.\]
We let $\mathcal{B}_k=\cup_{t=1}^{T_k}\mathcal{C}_t$ and finish the proof as before.

%%%%%%%%%%%%%%%%%%%%%%%%%%%%%%%%%%%%%%%%%%%%%%%%
\subsection{Classical means and dispersions} We remark that in Theorem \ref{th:p<d} we don't claim that $q_n$ and $r_n$ are generalized means of $b_n$ and $\widehat{b_n}$ respectively, for the definition of generalized means we refer the reader to \cite{JP}. However
for $p=2$ the construction above yields

\begin{corollary} 
For $d\ge 2$ there exists a basis $\{b_n\}_{n=1}^{\infty}$  for $L^2(\R^d)$  such that the sequences $\{\mu(b_n)\}_{n=1}^\infty$,
$\{\mu{\widehat(b_n)}\}_{n=1}^\infty$ and $\{\Delta(b_n)\Delta(\widehat{b_n})\}_{n=1}^\infty$ are bounded.
\end{corollary}

\begin{proof}
We repeat the construction described in the previous sections, but start with a function $\Psi^{(e)}=\Psi(x)+\Psi(-x)$. Then $\Psi^{(e)}_{j,s}$ has the same properties as $\Psi_{j,s}$ above with (\ref{eq:supp}) replaced with
\[
\supp(\Psi^{(e)}_{j,s})\subset\{x=(x_1,..., x_d): 2^{s-1}<x_m<2^s, m=1,...,d\}.\]
In addition we get $|\Psi^{(e)}_{j,s}(x)|=|\Psi^{(e)}_{j,s}(-x)|$ and
\begin{equation}
\label{eq:zeromean}
\mu(\Psi^{(e)}_{j,s})=\int_{\R^d}x|\Psi^{(e)}_{j,s}(x)|^2dx=0.
\end{equation}
Now for $d\ge 2$ and any $\beta_l$ in (\ref{eq:beta}) we have
\[
|\mu(\widehat{\beta_l})|\le\int_{\R^d}|\xi||\widehat{\beta_l}(\xi)|^2d\xi\le
\int_{\R^d}(2+|\xi-2^{-s}j_l|^2)|\widehat{\beta_l}(\xi)|^2d\xi\le C_1. 
\]
Clearly,
\[
\Delta(\beta_l)\Delta(\widehat{\beta_l})\le\int_{\R^d}|x|^2|\beta_l(x)|^2dx\int_{\R^d}|\xi-2^{-s}j|^2|\widehat{\beta_l}(\xi)|^2\le C_2.\]
Finally, to estimate $\mu(\beta_l)$ we use (\ref{eq:zeromean})   
\[
|\mu(\beta_l)|=|\mu(\beta_l)-\mu(\gamma_l\Psi^{(e)}_{l,s})|\le\]
\[\frac{\theta^2}{J_s}\int_{\R^d}|x||f(x)|^2dx+
\int_{\R^d}|x|\left||\sum_{n=1}^{l-1}\sigma_n\Psi^{(e)}_{n,s}(x)+\gamma_l\Psi^{(e)}_{l,s}(x)|^2-|\gamma_l\Psi^{(e)}_{l,s}(x)|^2\right|dx.\]
The first term is bounden by $J_s^{-1} I_2(f)\le 1$. (Here $g_t$ should replace $f$ for the case $d=2$.) We estimate the second term by
\[
(1+2^s)\int_{\R^d}|x||\sum_{n=1}^{l-1}\sigma_n\Psi^{(e)}_{n,s}(x)|^2dx+
2^{-s}\int_{\R^d}|x||\gamma_l\Psi^{(e)}_{l,s}(x)|^2dx\le\]
\[
(1+2^{s})2^{s}\int_{[-1,1]^d}|y||\sum_{n=1}^{l-1}\sigma_ne^{2\pi ij(n,s)\cdot y}\, \Psi^{(e)}(y)|^2dy
+\int_{[-1,1]^d}|y||\Psi^{(e)}(y)|^2dy\le 
\]
\[
2^{2s+1}\sqrt{d}\sum_{n=1}^{l-1}|\sigma_n|^2+\sqrt{d}\le \sqrt{d}(1+2^{2s+1}J_s^{-1})<C_3.\]
Thus $\mu(\beta_l)$ are bounded.
\end{proof}

%%%%%%%%%%%%%%%%%%%%%%%%%%%%%%%%%%%%%%%%%%%%%%%%%%%%%%%%%

\subsection{Degenerate homogeneous weights} We now consider uncertainty inequalities with homogeneous weights in $\R^d$, obtained by B.~Demange, see \cite{D}. We describe a simple construction of a bases of functions which have uniformly bounded weighted norms with there Fourier transforms. The construction is once again based on the argument form \cite{Bou}. This time a simple version of the argument implies
 
\begin{lemma} 
Let $u,v$ be non-negative functions on $\R^d$. 
Suppose that there exists an orthonormal sequence $\{\phi_j\}_{j=1}^\infty$ in $L^2(\R^d)$ of compactly supported smooth functions such that
\[
\int_{\R^d}u(x)|\phi_j(x)|^2 dx\le C_1^2\quad  {\rm and}\quad \int_{\R^d}v(\xi)|\widehat{\phi_j}(\xi)|^2 d\xi\le C_2^2,\]
for each $j$.
We assume also that the supports of $\phi_j$ form a locally finite family of sets; i.e. each compact set intersects only finite number of these supports.

Suppose also that there is a sequence $\{f_k\}_{k=1}^\infty$ that is dense in the unit sphere of $L^2(\R^d)$, functions $f_k$ are smooth and have compact supports and 
\begin{equation}
\label{eq:dense}
\int_{\R^d}u(x)|f_k(x)|^2 dx<\infty \quad  {\rm and}\quad \int_{\R^d}v(\xi)|\widehat{f_k}(\xi)|^2 d\xi< \infty,
\end{equation}
for each $k$.

Then for each $\epsilon>0$ there is an orthonormal basis $\{b_n\}_{n=1}^\infty$ for $L^2(\R^d)$ such that
\[
\int_{\R^d}u(x)|b_n(x)|^2 dx\le \left(C_1+\epsilon\right)^2 \quad  {\rm and}\quad \int_{\R^d}v(\xi)|\widehat{b_n}(\xi)|^2 d\xi\le \left(C_2+\epsilon\right)^2,\]
for each $n$.
\end{lemma}

Using this lemma we prove the following.

\begin{theorem}
Let $d>1$ and $\alpha_1,...,\alpha_d>0$ then there exists an orthonormal basis for $L^2(\R^d)$ such that
\[
\sup_n\int_{\R^d}|x_1|^{\alpha_1}|x_2|^{\alpha_2}...|x_d|^{\alpha_d}|b_n(x)|^2dx+
\int_{\R^d}|\xi_1|^{\alpha_1}|\xi_2|^{\alpha_2}...|\xi_d|^{\alpha_d}|\widehat{b_n}(\xi)|^2d\xi<\infty.\]
 
\end{theorem}

\begin{proof}
First, there exists a dense sequence $\{f_k\}_{k=1}^\infty$ in the unit sphere of $L^2(\R^d)$ that consist of smooth functions with compact supports; it satisfies (\ref{eq:dense}) for \[u(x)=v(x)=|x|^{\alpha_1}...|x_d|^{\alpha_d}.\]
To apply the lemma above we need to construct an orthonormal sequence bounded with their Fourier transforms in $L^2(\R^d,u(x))$.
Let $\phi$ be a smooth function such that its support is contained in the unit cube $[-1,1]^d$ and $\|\phi\|_2=1$.   We define
\[
\phi_j(x)=2^{\frac{j}{2}\left(\frac{\alpha_1}{\alpha_2}-1\right)}\phi(2^{-j}x_1-3,2^{\frac{j\alpha_1}{\alpha_2}}x_2,x_3,...,x_d).
\]
Then clearly 
$\|\phi_j\|_2=1$ and the support of $\phi_j$ is contained in the set 
\[
E_j=\{x\in \R^d: x_1\in [2^{j+1},2^{j+2}], |x_2|\le 2^{-j\frac{\alpha_1}{\alpha_2}}, |x_s|\le 1, s=3,...,d\},\]
these sets are disjoint and hence the sequence $\{\phi_j\}$ is orthogonal.
Further,
\[
\int_{\R^d}u(x)|\phi_j(x)|^2dx=\int_{E_j}|x_1|^{\alpha_1}...|x_d|^{\alpha_d}|\phi_j(x)|^2dx
\le 2^{(j+2)\alpha_1}2^{-j\alpha_1}\|\phi_j\|_2^2\le 2^{2\alpha_1}.\]
The Fourier transforms satisfy
\[
\widehat{\phi_j}(\xi)=e^{-2^{j+1}3\pi i\xi_1}2^{\frac{j}{2}\left(1-\frac{\alpha_1}{\alpha_2}\right)}\widehat{\phi}(2^j\xi_1,2^{-j\frac{\alpha_1}{\alpha_2}}\xi_2,...,\xi_d).\]
Thus
\[
\int_{\R^d}u(\xi|)\widehat{\phi_j}(\xi)|^2d\xi=
\int_{\R^d}|\xi_1|^{\alpha_1}|\xi_2|^{\alpha_2}...|\xi_d|^{\alpha_d}|\widehat{\phi}(\xi)|^2d\xi=C.\]
The theorem follows.    
\end{proof}

\section*{Acknowlegments}

The author is grateful to Yurii Lyubarskii for many valuable discussions of the material of the article. This work was done while the author was visiting Department of Mathematics of University of California, Berkeley, and it is a pleasure to thank the Department for its hospitality.

\end{document}